\documentclass[12pt]{book}
\usepackage{epsf,nemacom}
\newcommand{\D}[2]{\frac{\partial #2}{\partial #1}}

\renewcommand{\vec}[1]{\mbox{\boldmath$#1$}}
\newcommand{\Ord}[1]{{\cal O}\left(#1\right)}

\newcommand{\vth}{\vartheta}
\newcommand{\dx}{\partial_x}
\newcommand{\E}[1]{\mbox{e\footnotesize #1}\,}

\newtheorem{remark}{Remark}[section]

\begin{document}

\paperstart
%
%
%
%
%
%
\begin{flushleft}
%
%
\title{  Approximate Models of Dynamic Thermoviscoelasticity Describing
Shape-Memory-Alloy Phase Transitions }
%
%
%
\author{R.V.N.~Melnik and A.J.~Roberts}{Department of Mathematics and Computing, 
University of Southern Queensland, Toowoomba, \\ QLD 4305, Australia.}
%
%
%
%
\end{flushleft}
%
%
\markboth{R.V.N.~Melnik, A.J.~Roberts}{Approximate Models For Shape-Memory 
Alloys}
%
%
\index{Melnik, R.V.N.}
\index{Roberts, A.J.}
\index{Approximate Models for Shape-Memory Alloys}
\index{Nonlocal models of thermoviscoelasticity}
\index{Shape-memory alloys}
\index{Landau-Devonshire model}
\index{Centre manifold theory}
\index{Cattaneo-Vernotte equation}
\index{Phase transitions}
%
%

\section{Introduction}
In this paper we consider two models for the approximate description 
of thermomechanical behaviour  of viscoelastic materials.    
Accounting for thermal fields in such a description is important 
for all viscoelastic materials ranging from viscous fluids to 
elastic solids. The viscoelastic behaviour typically combines 
viscous and elastic properties and the relative proportion of this 
combination strongly depends on  thermal characteristics of the material. 
Moreover, with changing thermal conditions,  it is sometimes   difficult to 
decide whether 
a particular material is a solid or a fluid.
The key points in such 
decisions belong to the time of observation and  
 to the  choice of constitutive relations 
which couple stresses, deformation gradients, thermal fluxes 
and temperature.

Our analysis is based on  the nonlocal theory of 
continuum mechanics which considers constitutive variables defined at a point
as a function of 
their values over the whole spatial domain of interest rather 
than as a function at that point only \cite{Balta1977}. 
This  approach of  rational mechanics allows us to derive 
a general model that is suitable for the description of thermomechanical 
behaviour of materials under a wide range of temperature and loading 
patterns.   
In our models 
we allow for the dependency of stresses not only on 
the deformation gradient and temperature but also on the rates of their 
changes. Such  considerations put us closer to real situations where 
the time-dependent coupling between temperature and stresses 
have to include the velocity of the deformation gradient and the speed 
of thermal propagation.       
Another novelty of our paper is the  accounting for finite speeds of 
thermal disturbances. We define the constitutive relationship for thermal 
fluxes  using the Cattaneo-Vernotte equation which includes  the 
classical Fourier law  as a special limiting case 
(in the limit of zero relaxation time for heat fluxes). 
In particular this 
approach is critical  in modelling short 
transient states   in  low   
temperature regimes. 

During recent years a number of papers were devoted to 
the development of mathematical theory of thermomechanical phase 
transitions (see \cite{Niezgodka1988,Sprekels1989,Hoffmann1990,Hoffmann1995,Anderssen1998}
and references therein). 
The majority of those papers dealt with 
important theoretical issues of models such as well-posedness and
the global  asymptotic behaviour of
 solutions.  However, only a few papers have been devoted to the 
description of computational results using those models (see, for example,
\cite{Niezgodka1991,Klein1995} and references therein).
Almost all developed models  take into account neither 
the rate of  thermal disturbances nor the relaxation time of 
thermal fluxes.
However, the importance of these issues are well known in dynamic 
hyperbolic thermoelasticity where mathematical procedures and computational 
techniques have a longer history  compared to that in 
thermoviscoelasticity \cite{Melnik1997,Racke1997}.     
   
In dealing with the three main physical quantities of continuum mechanics 
(stresses, deformation 
gradients and displacements) it is important to take into 
account their coupling 
to the thermal field. This allows us to construct efficient 
mathematical models for the description of complicated phenomena,
 such as 
hysteresis,
which  are becoming increasingly important in a wide range of applications. 
In this paper we apply the developed models to the description of 
shape memory alloy effects in a large bar. 
It is  well-known  that for many types of shape memory 
 materials the  dependency
of stresses on the deformation gradient  
upon loading and 
unloading is significantly different.
Applying a large load at a low temperature, we may get a residual 
deformation gradient, which typically  vanishes upon heating. 
The restoring of the original shape is  referred to as the 
shape memory effect. This effect is discussed with two numerical 
examples.

The rest of the  paper is organized as follows. 
\begin{itemize}
\item Section 2 provides the reader with basic preliminaries and notation.
\item  The general formulation of the model  is given in Section 3.
In this section we specify the model for internal energy and derive 
restrictions on the model imposed by the second law of thermodynamics.
\item
In Section 4 we incorporate the Cattaneo-Vernotte equation for heat conduction
into our model.
\item
Section 5 deals with the Landau-Devonshire model for the free energy function.
The constitutive relation connecting stresses and the deformation gradient
is also discussed in this section.
\item
In Section 6 we consider a one-dimensional model of thermoviscoelasticity 
and discuss the consequences of non-convexity of free energy
function.
\item
Some numerical results are presented and discussed  in Section 7.
\item
In Section 8 we use centre manifold theory to derive
an approximate mathematical model for the description
of thermomechanical behaviour of viscoelastic materials.
\end{itemize}

\section{Preliminaries and Notation}
Assume that an object of interest (a solid, fluid, gas or plasma) 
 occupies the volume $V$  in a 
fixed reference spatial configuration $\Omega$ at a certain time $t_0$. 
This object in its spatio-temporal configuration will be referred to
by the generic name ``system''.
 We aim to develop an  efficient mathematical 
description of the dynamic thermomechanical behaviour of the system.

Let ${\bf x}=(x_1, x_2, x_3)$ be material (Lagrangian) 
coordinates of a material point of the system in the configuration 
$\Omega$ at time $t_0$. Then the dynamics of the system  is
 determined by the spatial displacements ${\bf u}=(u_1, u_2, u_3)$ 
of such material points 
as a function of the reference position, ${\bf x}$, 
and the time of interest, $t$. 
The partial derivative of displacement with respect to ${\bf x}$
is identified with the symmetric strain tensor   
\begin{eqnarray}
\mbox{\boldmath $\epsilon$} = {\rm sym} \left [ \frac{\partial {\bf u}
({\bf x}, t)}{\partial {\bf x}} \right] 
\quad \mbox{or} \quad 
\epsilon_{ij}({\bf x}, t) = \frac{1}{2} \left 
[\frac{\partial u_i({\bf x}, t)}{\partial x_j} +  
\frac{\partial u_j({\bf x}, t)}{\partial x_i} \right]
\, , \quad i, j=1,2,3 \,, 
\label{2eq1}
\end{eqnarray}
and the time derivative of the function ${\bf u}$ is identified with
 the velocity of 
the system
\begin{eqnarray}
{\bf v} = \frac{\partial {\bf u}}{\partial t} \quad
\mbox{or} \quad v_i({\bf x}, t) = 
\frac{\partial u_i({\bf x}, t)}{\partial t} \, , \quad  i=1,2,3 \,.
\label{2eq2}
\end{eqnarray}
In (\ref{2eq1}) we  require that  ${\rm det} 
(I + \mbox{\boldmath $\epsilon$}) >0$ which 
precludes a possibility of compression of the matter to zero 
and guarantees the local invertibility of ${\bf x} + {\bf u}({\bf x}, t)$
 \cite{Renardy1987}. 
Since the time derivatives   are understood in the 
Lagrangian sense,  ${\bf x}$ is kept fixed in (\ref{2eq2}).

\section{ Nonlocal Models of Thermoviscoelasticity}
The equation of motion requires information on forces acting per unit area
of the matter and, hence, in a  natural way, involves 
 the concept of stresses. The stress is not a mere function of
of the deformation gradient, as it is often assumed. It 
 also depends on temperature of the matter, its rate of change in time
and the rate of change of 
deformation gradient $\mbox{\boldmath $\epsilon$}$. 
Let $\rho_0({\bf x}, t_0)>0$ be the density of the matter (the mass 
per unit volume)
in the reference configuration $\Omega$ at time $t_0$ and $\rho({\bf x}, t)$
be the density of the matter at time $t$ where $t-t_0$ is sufficiently small.
Then, in the Lagrangian system of coordinates 
$({\bf x}, t)$, the equation for balance of mass is written 
in the form \cite{Renardy1987}
\begin{eqnarray} 
\rho({\bf x}, t) {\rm det} (I+ \vec\epsilon({\bf x}, t)) =
\rho_0({\bf x}, t_0) \,. 
\label{2eq11}
\end{eqnarray}
The  equation of motion has the following form
\begin{eqnarray}
\rho \frac{\partial^2 {\bf u}}{\partial  t^2} = \nabla_{{\bf x}} \cdot 
{\bf s} + {\bf F} \quad \mbox{with} \quad 
{\bf F} = \rho({\bf f} + \hat{{\bf f}}) - {\hat{\rho}} {\bf v} \,,
\label{2eq12}
\end{eqnarray}
where ${\bf f}$ is a given body  force per unit mass,
$ \hat{\rho}$ and $\hat{{\bf f}}$ are nonlocal mass and force residuals
respectively,  and $\bf{s}$ is the 
stress tensor. 

In  Lagrangian coordinates the equation for energy balance
has the form 
\begin{eqnarray}
\rho \frac{\partial}{\partial t} \left (e + \frac{ {\bf v}^2 }{2} \right )
- \nabla_{{\bf x}} \cdot( {\bf s} \cdot {\bf v} ) + \nabla \cdot {\bf q} =
\rho \left ( h + \hat{h} + {\bf f} \cdot {\bf v} - 
\frac{ \hat{\rho} }{\rho} \left ( e +
 \frac{ {\bf v}^2 }{2} \right ) \right ),
\label{2eq16}
\end{eqnarray}
where $e$ is the specific internal energy of the system, 
${\bf v}^2 = {\bf v} \cdot 
{\bf v}$, $h$ is the heat source density, $\hat{h}$ is the nonlocal energy 
residual (see \cite{Balta1977} for conditions on localised residuals) 
and ${\bf q}$ is the heat flux.
The scalar multiplication of (\ref{2eq12}) by ${\bf v}$ gives
\begin{eqnarray}
\rho \frac{ \partial {\bf v}^2/2 }{\partial t} - {\bf v} \cdot (\nabla \cdot 
{\bf s}) = 
({\bf F}, {\bf v}) \equiv
\rho ( {\bf f} + \hat{{\bf f}})  \cdot {\bf v}  - { \hat{\rho} } {\bf v}^2.
\label{2eq17}
\end{eqnarray}
Taking into account normalisation, from (\ref{2eq16}) and (\ref{2eq17}) we get
\begin{eqnarray}
\rho \frac{\partial e}{\partial t} - {\bf s}^T : (\nabla {\bf v}) + \nabla \cdot
 {\bf q}
= g \, ,
\label{2eq18}
\end{eqnarray}
where $\displaystyle {\bf a}^T : {\bf b} = \sum_{i, j =1}^3 a_{ij} b_{ij}$ is the 
standard notation for the rank~2 tensors 
${\bf a}$ and ${\bf b}$  and
\begin{eqnarray}
g= \rho (h + \hat{h}) - \rho \hat{{\bf f}} \cdot {\bf v} - 
{\hat{\rho}}
\left ( e - \frac{ {\bf v}^2 }{2} \right ).
\label{2eq19}
\end{eqnarray}   
The right-hand sides of equations (\ref{2eq12}) and (\ref{2eq18}) 
incorporate into the model 
nonlocal and dissipative effects of thermomechanical waves.
As we shall see in the next section, under appropriate constitutive 
relations it is also possible to  allow
for a relaxation time for acceleration of the motion in response to 
applied gradients such as the deformation gradient and the 
temperature gradient.


We assume that  there exists  a one-to-one entropy function  
of the system state. We denote the density of such a function
 by $\eta$, and then the second law of thermodynamics is
\begin{eqnarray}
\frac{\partial \eta}{\partial t} - \nabla \cdot {\bf r} \geq \xi + 
\hat{\xi} - 
\frac{ \hat{\rho} }{\rho} \, ,
\label{2eq20}
\end{eqnarray}
where $\xi$ is the entropy source density, ${\bf r}$ is the entropy flux density
and $\hat{\xi}$ is the nonlocal entropy residual.

The system of equations (\ref{2eq12}), (\ref{2eq18}) combined with 
inequality (\ref{2eq20}) 
provides the general  mathematical model for the 
description of thermomechanical behaviour of dynamic systems.
The macroscopic modelling of such systems starts from the choice of 
constitutive relationships. We 
assume the existence of a functional $\Psi$ invariant under a 
time shift and 
chose this functional in the form of the 
Helmholtz free energy 
\begin{eqnarray}
\Psi = e - \theta \eta \, ,
\label{2eq21}
\end{eqnarray}
 where $\theta$ is the temperature of the system ($\theta >0,$
$\inf_{({\bf x},t)} \theta =0$).  We also assume 
specific forms for the entropy flux and the entropy source density as 
\begin{eqnarray}
{\bf r} = {\bf q}/ \theta \, , \quad \xi = h/ \theta \,.
\label{2eq22}
\end{eqnarray}
Using  (\ref{2eq21})  in (\ref{2eq18}) and taking into account
that 
\begin{eqnarray}
\nabla \cdot  {\bf q} = \theta \nabla \cdot ({\bf q}/ \theta) 
+ ( {\bf q} \cdot \nabla 
\theta)/\theta \, ,
\label{2eq23}
\end{eqnarray}
from (\ref{2eq20}) and (\ref{2eq22}) 
we get the nonlocal formulation of the Clausius-Duhem inequality
\begin{eqnarray}
 - \frac{ \hat{\rho} }{\rho} \left ( \Psi - \frac{ {\bf v}^2 }{2} \right )
- \left ( \frac{ \partial \Psi}{\partial t} + \eta 
\frac{ \partial \theta }{\partial t}  \right ) + {\bf s}^T : \nabla {\bf v} - 
\hat {{\bf f}} \cdot {\bf v} - \frac{ {\bf q} \cdot \nabla \theta }{\theta}
- (\theta \hat{\xi} - \hat{h} ) \geq 0 \,.
\label{2eq24}
\end{eqnarray}
The latter inequality together with requirements on localisation residuals
(see \cite{Balta1977} for details) impose
restrictions on the choice of 
nonlocal residuals
and the functions $\eta$, ${\bf s}$ and ${\bf q}$.
We assume that the entropy density is  given in the form 
\begin{eqnarray}
\eta= - \frac{\partial \Psi}{\partial \theta} \,.
\label{2eq25}
\end{eqnarray}
Finally, we have to specify
the constitutive relationships that couple  stresses,  
deformation gradients, temperature and heat fluxes
\begin{eqnarray}
\Phi_1({\bf s}, \mbox{\boldmath $\epsilon$})=0 \,, \quad
\Phi_2( {\bf q},  \theta)=0 \, ,
\label{2eq26}
\end{eqnarray}
where it is implicitly assumed that 
these relations may involve spatial and temporal 
derivatives of the functions.
In Section 4 and 5 we specify particular forms for $\Phi_1$ and $\Phi_2$.

\section{The Cattaneo-Vernotte Model for  Heat Conduction}
The choice of the function $\Phi_2$ in (\ref{2eq26}) is made using  
the  Cattaneo-Vernotte model
\begin{eqnarray}
{\bf q} + \tau_0 \frac{\partial  {\bf q} }{\partial t} = 
- k(\theta, \mbox{\boldmath $\epsilon$}) \nabla \theta \, ,
\label{3eq1}
\end{eqnarray}
where $\tau_0$ is the dimensionless thermal relaxation time
and  $k(\theta, \vec\epsilon)$ is the thermal conductivity of the 
material (typically $k=1 +\tilde{\beta} \theta$ with the given 
dimensionless coefficient $\tilde{\beta}$).
Such a choice is made in order to account for the finite speeds of 
thermal wave propagation and thermally induced stress wave propagation
coupled to the deformation gradient \cite{Glass1994,Muller1993}. 
In order to incorporate  equation (\ref{3eq1}) into the general model 
of thermoviscoelasticity  we use a consequence of (\ref{2eq18})
\begin{eqnarray}
\rho \tau_0 \frac{ \partial^2 e }{\partial t^2} - \tau_0 
\frac{ \partial }{\partial t} [ {\bf s}^T : ( \nabla {\bf v})]
+ \tau_0 \nabla \cdot \left ( \frac{\partial {\bf q}}{\partial t} \right )
= \tau_0 
\frac{ \partial g}{\partial t} \, .
\label{3eq2}
\end{eqnarray}
On the other hand, from (\ref{3eq1}) we get
\begin{eqnarray}
\nabla \cdot {\bf q} + \tau_0 \nabla \cdot 
\left ( \frac{\partial {\bf q}}{\partial t} \right )
=- \nabla \cdot (k \nabla \theta) \, .
\label{3eq3}
\end{eqnarray}
Then from (\ref{2eq18}), (\ref{3eq2}), (\ref{3eq3}) we obtain
the energy balance equation in the form
\begin{eqnarray}
\rho \frac{\partial e}{\partial t} + \rho \tau_0 \frac{ \partial^2 e }
{\partial t^2} 
- {\bf s}^T :( \nabla {\bf v}) - 
\tau_0 \frac{\partial}{\partial t}[{\bf s}^T :( \nabla {\bf v})]
- \nabla \cdot (k \nabla \theta) = 
G \, ,
\label{3eq4}
\end{eqnarray}
where
\begin{eqnarray}
G=g + \tau_0 \frac{ \partial g}{\partial t} \, .
\label{3eq5}
\end{eqnarray}
During recent years, the interest in such a hyperbolic approach in 
the analysis of materials with memory has  increased 
 \cite{Colli1995}.

\section{The Landau-Devonshire Model for the Helmholtz Free Energy 
and the Stress-Strain Relation}
We start from the consideration of the one-dimensional case assuming
 the following  approximation 
 for the  
free energy of the system
\begin{eqnarray}
\Psi(\theta, \epsilon) = \psi_0(\theta) + \psi_1(\theta) 
\psi_2(\epsilon) + \psi_3(\epsilon)
\label{4eq1}
\end{eqnarray}
where $\psi_0(\theta)$  models   
thermal field contributions,  $\psi_1(\theta) 
\psi_2(\epsilon)$  models  shape-memory contributions and 
 $\psi_3(\epsilon)$  models  mechanical field contributions.
These models are chosen in the following forms
\begin{eqnarray}
\left [
\begin{array}{l}
\psi_0(\theta) = \alpha_0 - \alpha_1 \theta \ln \theta \, , \quad
\psi_1(\theta)= \frac{1}{2} \alpha_2 \theta \, , \quad
 \psi_2(\epsilon) = \epsilon^2 \, , 
\\[10pt]
\displaystyle
\psi_3(\epsilon) = -\frac{1}{2} \alpha_2 \theta_1 
\epsilon^2 
-\frac{1}{4} \alpha_4  \epsilon^4+ 
\frac{1}{6} \alpha_6  \epsilon^6 \, ,
\end{array}
\right.
\label{4eq2}
\end{eqnarray}
where  all $\alpha_i$ and $\theta_1$ are positive constants.
The model (\ref{4eq1})--(\ref{4eq2}), known as the Landau-Devonshire 
model for the Helmholtz 
free energy,  
 covers a number of important practical cases. 
However, it belongs to the class of models  which 
is difficult to investigate compared to the 
Landau-Devonshire-Ginzburg model. In the latter case  
an additional ``smoothing'' term in (\ref{4eq1}), known as the Ginsburg term
$\gamma u_{xxxx}$, allows us to obtain 
a bound of  the deformation 
gradient (strain) using a well established technique \cite{Chen1994}.

\begin{remark}
A number of important characteristics of phase transformations (such as 
the size of hysteresis) may depend 
on the contributions of the interfacial energies. These contributions 
are often modeled with the Ginsburg correction term.
However, the Ginsburg coefficient can only be determined in approximate order
\cite{Sprekels1990} and in the general case this coefficient 
may not be temperature-independent. 
Another way to account for the contributions of interfacial energies is 
to to take the free energy in the form \cite{Bornert1990,Moyne1997}
\begin{eqnarray}
\Psi =(1-z) \tilde{\psi}_1 (\epsilon, \theta)+z \tilde{\psi}_2 
(\epsilon, \theta) +z(1-z) \tilde{\psi}_3 \, , 
\end{eqnarray}
where $z$ is the volume fraction of martensite (i.e. the product phase),
$(1-z)$ is the volume fraction of  austenite (i.e. the parent phase),
$\tilde{\psi}_1$, $\tilde{\psi}_2$ are the free energies of austenite and martensite respectively and  $\tilde{\psi}_3$ is the contribution from the interaction effect between austenite and martensite.
We will not pursue these ideas  in this paper.
\end{remark}
\begin{remark}
Some authors include a linear term $\alpha^0 \theta$ into $\psi_0(\theta)$.
This term has  no bearing on the final model and changes only the value 
of the coefficient of  $\theta$ in the internal energy representation 
(see formula (\ref{5eq2})), and thus is omitted. 
\end{remark}  

In the general case for the choice of the function 
$\Phi_1$ in (\ref{2eq26}) we 
allow the dependency of the stress on the rate of  
temperature and the deformation gradient
\begin{eqnarray}
 s = \rho \left [p(\theta, \epsilon) + 
\lambda \left ( \frac{\partial \theta}{\partial t} \, , \quad   
\frac{\partial \epsilon}{\partial t} \right ) \right ] ,  
\label{4eq3}
\end{eqnarray}
where 
\begin{eqnarray}
p(\theta, \epsilon) = \frac{\partial \Psi}
{\partial \epsilon} \, , \quad
\lambda \left ( \frac{\partial \theta}{\partial t},  
\frac{\partial \epsilon}{\partial t} \right ) = 
\tilde{\mu}(\theta) \frac{ \partial \epsilon}{\partial t} + 
\tilde{\nu}(\epsilon) 
\frac{\partial \theta}{\partial t} \, .
\label{4eq4}
\end{eqnarray}
It is straightforward to deduce
\begin{eqnarray}
p(\theta, \epsilon) = \alpha_2 
\theta \epsilon 
+ \frac{ \partial \psi_3 (\epsilon) }{\partial 
\epsilon }=
\left [\alpha_2  \epsilon (\theta - \theta_1) - 
\alpha_4 \epsilon^3
+ \alpha_6 \epsilon^5 \right ].
\label{4eq5}
\end{eqnarray}

\section{One-Dimensional Hyperbolic Approximation of Shape-Memory-Alloy 
Dynamics} 
Using the model (\ref{4eq1}) and (\ref{4eq2}), from (\ref{2eq25}) we get
\begin{eqnarray}
\eta = \alpha_1 (1+ \ln \theta) - \frac{1}{2} \alpha_2 \epsilon^2 \, .
\label{5eq1}
\end{eqnarray}
This enables us to find the internal energy of the system as
a sum of thermal and mechanical fields contributions
\begin{eqnarray}
e= \alpha_0 + \alpha_1 \theta - 
\frac{1}{2} \alpha_2 \theta_1 \epsilon^2 
- \frac{1}{4} \alpha_4 \epsilon^4+ 
\frac{1}{6} \alpha_6 \epsilon^6 =
\alpha_0 + \alpha_1 \theta + \psi_3 (\epsilon) \, .
\label{5eq2}
\end{eqnarray}
The substitution of  (\ref{5eq2}) into (\ref{3eq4})  leads 
to the final form of the energy balance equation.
In particular, assuming  symmetry of the deformation gradient tensor,
we get 
\begin{eqnarray}
\rho \alpha_1 \left [ \frac{\partial \theta}{\partial t} + \tau_0 
\frac{\partial^2 \theta}{\partial t^2} \right ] +
A(\epsilon, \theta) - \nabla \cdot (k \nabla \theta) =G \, ,
\label{5eq3}
\end{eqnarray}
where the meaning of $A$ is
\begin{eqnarray}
\displaystyle
& &
A(\epsilon, \theta) =
- \rho \alpha_2 \left \{ \theta \epsilon 
 \frac{ \partial \epsilon }
{\partial t} + \tau_0 \frac{\partial }{\partial t}
\left [ \theta \epsilon 
 \frac{ \partial \epsilon }
{\partial t} \right] \right \} -
\rho \tilde{\mu}(\theta) \left \{ \left ( \frac{ \partial \epsilon }
{\partial t} \right )^2  + 
\right.
\nonumber
\\[10pt]
\displaystyle
& &
\left.
\tau_0 \frac{\partial }{\partial t} \left [
\left ( \frac{ \partial \epsilon }
{\partial t} \right )^2  \right ] \right \} -
\rho \frac{\partial \theta}{\partial t} 
\left \{ \tilde{\nu}( \epsilon) 
\frac{ \partial \epsilon }
{\partial t} + \tau_0 \frac{\partial }{\partial t} \left [
\tilde{\nu}( \epsilon )
\frac{ \partial \epsilon }
{\partial t} \right ] \right \}.
\label{5eq4}
\end{eqnarray}
Equation (\ref{5eq3}) is solved together with the equation of motion 
(\ref{2eq12}) with respect to $( u, \theta)$:
\begin{eqnarray}
\left \{
\begin{array}{l}
\displaystyle
C_v \left [ \frac{\partial \theta}{\partial t} + \tau_0 
\frac{\partial^2 \theta}{\partial t^2} \right ]
- k_1 \left [
\theta \frac{\partial u}{\partial x} \frac{\partial^2 u}
{\partial t \partial x} + \tau_0 \frac{\partial}{\partial t}
\left ( \theta \frac{\partial u}{\partial x} \frac{\partial^2 u}
{\partial t \partial x} \right ) \right ] -
\mu \left [ \left ( \frac{\partial^2 u}{\partial t \partial x} 
\right )^2 +
\right.
\\[10pt]
\displaystyle
\left.
\tau_0 \frac{\partial}{\partial t} \left ( \frac{\partial^2 u}
{\partial t \partial x} \right )^2 \right ] -
 \nu \left [ \frac{\partial \theta}{\partial t}
\frac{\partial^2 u}
{\partial t \partial x} + \tau_0 
\frac{\partial }{\partial t} \left ( \frac{\partial \theta}{\partial t}
\frac{\partial^2 u} {\partial t \partial x} \right ) \right ] - 
\frac{\partial}{\partial x} \left ( k \frac{\partial \theta}{\partial x}
\right ) = G \, ,
\\[10pt]
\displaystyle
\rho \frac{\partial^2 u}{\partial t^2} - 
\frac{\partial}{\partial x} 
\left [ k_1  \frac{\partial u} { \partial x} (\theta - \theta_1) 
- k_2  \left (\frac{\partial u} { \partial x} \right )^3 +
k_3  \left ( \frac{\partial u} { \partial x} \right )^5 \right ]
- \mu \frac{\partial^3 u}{\partial x^2 \partial t} - 
\nu \frac{\partial^2 \theta}{\partial x \partial t} = F \, ,
\end{array}
\right.
\label{5eq5}
\end{eqnarray}
where $C_v=\rho \alpha_1$, $k_1 = \rho \alpha_2$, $k_2 = \rho \alpha_4$,
$k_3= \rho \alpha_6$, $ \mu= \rho \tilde{\mu}$, $\nu = \rho \tilde{\nu}$.

The initial conditions for the model (\ref{5eq5}) 
are chosen in the form
\begin{eqnarray}
u(x, 0) = u^0(x) \, , \quad \frac{\partial u}{\partial t} (x,0) = u^1(x) \, ; 
\quad
\theta(x, 0) = \theta^0(x) \, , \quad \frac{\partial \theta}{\partial t} 
(x,0) = 
\theta^1(x) \, ,
\label{5eq6}
\end{eqnarray}
for given functions $u^0$, $u^1$, $\theta^0$, $\theta^1$. 
There are several distinct choices for  boundary conditions to be used in our computational experiments.
Mechanical boundary conditions are taken in one of the following forms
($L$ is the length of the structure): 
\begin{itemize}
\item  ``stress-free'' boundary conditions: $s(0, t)= s(L, t) =0$;
\item ``pinned end'' boundary conditions:
$u(0, t) = u(L, t) =0$;
\item or mixed mechanical boundary condition: $s(0, t)=0, \; \; u(L, t)=0$.
\end{itemize}
When displacements are given on boundaries, {\em a priori} bounds on 
strains are generally unknown 
which complicates the mathematical analysis of the problem. Computational 
results presented in Section 7 deal with this case. 
Thermal boundary conditions are chosen in one of the following form
\begin{itemize}
\item ``thermal insulation'' boundary conditions:
$ q(0, t)= q(L, t) =0 \, ,$
which reduce to  
$\displaystyle \frac{\partial \theta}{\partial x} (0, t) =
\frac{\partial \theta}{\partial x} (L, t) = 0$
for the Fourier law (see (\ref{6eq5}));
\item ``controlled flux''  boundary conditions: 
$\displaystyle \frac{\partial \theta}{\partial x} 
(0, t) =0 \, , \quad
- k \frac{\partial \theta}{\partial x} (L, t) = \beta[\theta - 
\theta^0(t)]$;
\item or fixed temperature (``uncontrolled energy flow'') boundary conditions:
$\theta(0, t)= \theta(L, t)=0$.
\end{itemize}
In the last case  additional assumptions are needed.
By using the Leray-Schauder principle we have  analysed 
 the Cauchy problem for 
nonlinear hyperbolic model of thermoviscoelasticity  
(\ref{5eq5}). 
Our procedure makes use of  
the Lumer-Phillips theorem and the technique developed  in 
\cite{Chen1994}. We shall address details of this 
technique  elsewhere.

Our final remark in this section goes to 
the choice of the function $\Psi$ in the form (\ref{4eq1}) and
 (\ref{4eq2})
that
  brings major 
difficulties in the investigation of the model (\ref{5eq5}).
Strictly speaking, the free energy
 function strongly depends upon
 the statistics of the phenomenon and has to be derived from  a statistical model. 
Since van der Waals work on statistical mechanics it is a common practice
to choose  this function as a  non-convex function of 
$\mbox{\boldmath $\epsilon$}$
\cite{Huo1994}.
When dealing with shape memory alloys, minima of this function are known to correspond different phases of the 
material. For example, in the case of three minima, we  expect one 
austenitic and two martensitic phase (see, for example, 
 \cite{Falk1980,Niezgodka1988,Sprekels1990,Hoffmann1990}).  
Temperature plays a crucial role in the phase transition. Depending 
on the  value of temperature, the material may alternate between a single 
thermodynamically unstable nonmonotone branch and multiple unstable branches. 
The character of this instability depends not only on the deformation 
gradient and temperature, but also on the rates of their changes.

\section{Computational Experiments}
In this section  we 
present some numerical results on the thermal and mechanical 
control  of a rod ($L=1 {\rm cm}$) 
with a shape-memory-alloy
core. The parameters of the Cu-based core are taken as follows 
\begin{eqnarray*}
& &
\displaystyle
k=1.9 \times 10^{-2} {\rm cm g/(ms^3 K)}, \quad \rho =11.1 {\rm g/cm^3},
\quad C_v =29 {\rm g/(ms^2 cm K)}, \quad \theta_1 =208 {\rm K} \, , 
\\[10pt]
& &
\displaystyle
k_1 =480 {\rm g /(ms^2 cm K)}, \quad 
k_2 = 6 \times 10^6 {\rm g /(ms^2 cm K)}, \quad
k_3= 4.5 \times 10^8 {\rm g /(ms^2 cm K)} \, . 
\end{eqnarray*}
We use  model (\ref{5eq5}) with $\tau_0=0= \mu= \nu=0,$  
 initial conditions (\ref{5eq6}) and 
``pinned end \& controlled flux'' boundary conditions.
This model was straightforwardly reduced to a differential-algebraic system 
in $\vec\pi=(u, v, \theta)^T$ and stress $s$ using second-order accurate 
spatial differences on staggered grids:
\begin{eqnarray}
D \frac{\partial \vec\pi}{\partial t} = {\bf f} \, , 
\quad s=k_1(\theta - \theta_1) 
\frac{\partial u}{\partial x} - k_2 
\left(\frac{\partial u}{\partial x} \right )^3 + k_3 
\left(\frac{\partial u}{\partial x} \right )^5 \, ,
\end{eqnarray}
where $D$ is the diagonal matrix with ${\rm diag}(D) = (1, \rho, C_V)$,
${\bf f} = (f_1, f_2, f_3)^T$ and
\begin{eqnarray}
f_1 =v \, , \quad f_2 = \frac{\partial s}{\partial x} + 
 F \, , \quad 
f_3 = k \frac{\partial^2 \theta}{\partial x^2} + k_1 \theta 
\frac{\partial u}{\partial x} \frac{\partial v}{\partial x} +G \, .
\label{8eq2}
\end{eqnarray}
The developed code is robust and much simpler compared to 
computational procedures previously reported in the literature 
for shape-memory alloys \cite{Niezgodka1991,Klein1995,Hoffmann1995}.

{\bf Experiment 1 (thermal control of phase transformations)}. 
In this experiment we set uniform forcing  $F=500 {\rm g/(ms^2 cm^2)}$
and vary heating conditions 
given by $G = 375 \pi \sin^3 (t \pi/6) {\rm g/ (ms^3 cm)}$.
We assume that the initial displacements are given in the form
\begin{eqnarray}
u^0(x) = 
\left \{
\begin{array}{ll}
-0.11809x, & 0 \leq x \leq 1/6, \\
0.11809(x-1/3), & 1/6 \leq x \leq 1/2, \\
0.11809(2/3-x), & 1/2 \leq x \leq 5/6 \\
0.11809(x-1), & 5/6 \leq x \leq  1
\end{array}
\right.
\end{eqnarray}
and take the initial temperature 
as $\theta^0=200 {\rm K}$. Figure 1 (obtained with time step 
$7 \times 10^{-4} {\rm ms}$ and 
 space step $1/24 {\rm cm}$) demonstrates the transformation of  
 $2 M^+ + 2 M^-$ martensites into 
an austenite (visible in the region of zero strain and 
displacements with superposed elastic vibrations as seen most clearly in the velocity field)
after sufficient temperature 
has reached. Then upon cooling we  observe 
a first order 
(martensitic) transition from the high temperature phase (austenite)
to the low temperature phase (martensite).
\begin{figure}[h]
\renewcommand{\epsfsize}[2]{0.8#1}
\centerline{\epsfbox{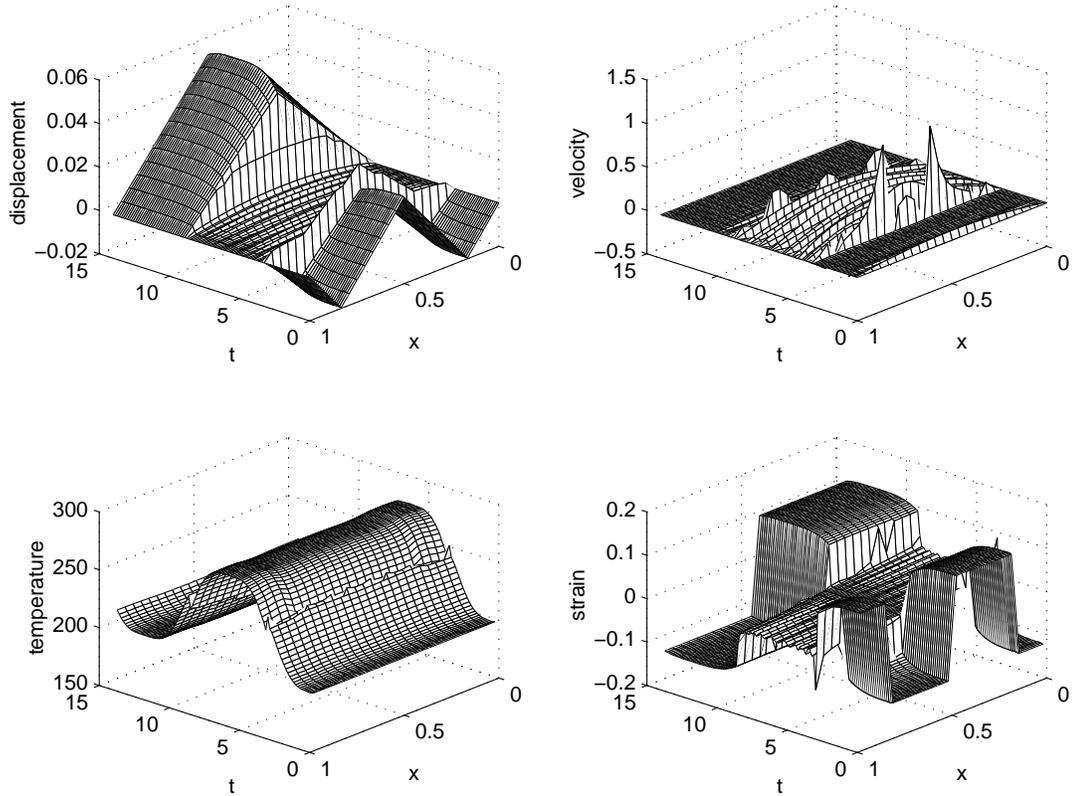}}
\caption{Thermally induced phase transformations.}
\end{figure}
Upon the return to the  low temperature regime 
the stable attractor with this applied thermomechanical forcing  
is not the original configuration but only  two distinct 
martensite phases. The transformation 
$[2 M^+ + 2 M^-] \rightarrow A$ is accompanied by a decrease in temperature    
whereas the transformation $A \rightarrow [M^+ + M^-]$ is accompanied by 
an increase in temperature.

{\bf Experiment 2 (mechanical control of phase transformations.} 
In this experiment we set $G=0$, but vary the mechanical loading according 
to  $F= 7000 \sin^3(\pi t/2)$ ${\rm g/(cm^2 ms^2)}$. 
Starting from the austenite configuration ($u^0=0$)
at intermediate temperature $\theta^0=255 {\rm K}$ we observe 
(see Figure 2 where the time step was  $8 \times 10^{-4} {\rm ms}$
 and the space step was $1/16 {\rm cm}$) 
the transformation $A \rightarrow [M^+ + M^-] \rightarrow A 
\rightarrow [M^- + M^+] \rightarrow A$. 
In this experiment  we observe the almost immediate transformations 
of austenite phases into 
two martensites 
  upon the increase/decrease in 
loading. 
Note the relatively large heating/cooling associated with the transition into/out of martensite phase.
A similar behaviour under different thermomechanical conditions 
was also observed in  \cite{Niezgodka1991,Klein1995}.
\begin{figure}[h]
\renewcommand{\epsfsize}[2]{0.8#1}
\centerline{\epsfbox{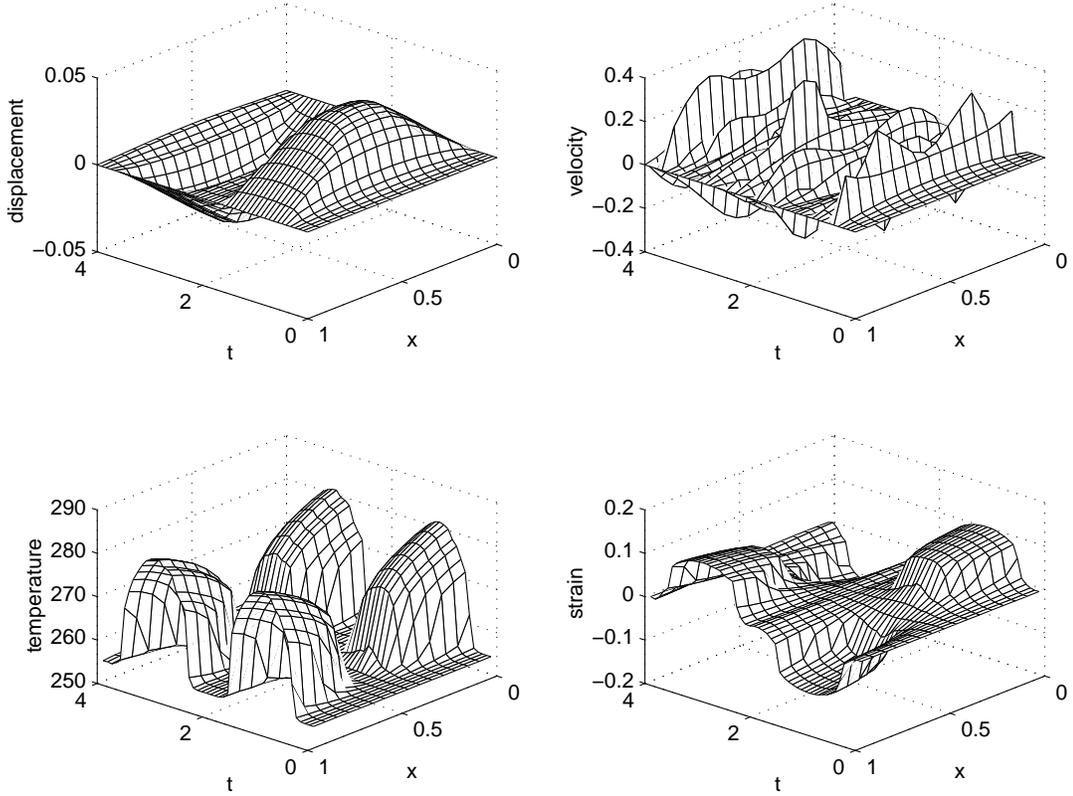}}
\caption{Mechanically induced phase transitions.}
\end{figure}
In our code we have also incorporated the Ginsburg term by
adding  $\gamma u_{xxxx}$ to $f_2$ in (\ref{8eq2}). 
With reported values of the Ginsburg coefficient ($\gamma \sim
10^{-10}-10^{-12}$) the Ginsburg term has a negligible effect 
on the thermomechanical behaviour of shape-memory alloys in the group 
of experiments described here. Accounting for interfacial energy 
contributions and the  influence of mechanical 
and thermal dissipations on the dynamics of memory material 
require further investigation.

\section{Construction of Approximate Models for Dynamic 
Thermoviscoelasticity Using Centre Manifold Theory} 
The model described in Section 6 will provide a good 
approximation of  thermomechanical behaviour of 
a large shape memory alloy bar (see applications in 
\cite{Besselink1997}) only in the case 
 the bar 
can be modelled by a thin rod with a shape memory alloy core.
As an alternative  to that model, in this section  we construct a new model 
which is derived  directly from the 3D model 
for shape memory alloy evolution (see (\ref{2eq12}), (\ref{5eq3})) 
using  centre manifold techniques (see, for example, \cite{Roberts1993}).

We assume that  the shear stress in equation (\ref{2eq12}) 
is determined by its three components,
the quasi-conservative component, ${\bf s}^q$,
the stress component due to mechanical dissipation,
${\bf s}^m,$
and the stress component due to thermal dissipations,
${\bf s}^t,$ (the latter is assumed to be negligible at this stage)
\begin{eqnarray}
{\bf s}= 
{\bf s}^q + {\bf s}^m + 
{\bf s}^t, \quad \mbox{with} \quad {\bf s}^q= \rho \frac{\partial \Psi}
{\partial \mbox{\boldmath $\epsilon$}},
\quad
{\bf s}^m = \rho \mu 
\frac{\partial \mbox{\boldmath $\epsilon$}}{\partial t}, \quad 
{\bf s}^t={\bf 0}.
\label{6eq3}
\end{eqnarray}

In the general case the heat flux is  determined as the 
solution of equation 
(\ref{3eq1}).  An approximation to this solution
is provided by the following generalised form
(see \cite{Niezgodka1988} and references therein)
\begin{eqnarray}
{\bf q} = -k \nabla \theta - \alpha \frac{\partial k \nabla \theta}
{\partial t}, \quad \alpha \geq 0, 
\label{6eq5}
\end{eqnarray}
 which we will use
with $\alpha=0$ when (\ref{6eq5}) turns into the classical Fourier 
law.

The internal energy function $e$ is defined from (\ref{2eq21}) by
\begin{eqnarray}
e = \Psi - \theta \frac{\partial \Psi}{\partial \theta}.
\label{6eq6}
\end{eqnarray}
In order to complete the formulation of the problem 
we  specify  a model for 
the free energy
function $\Psi$. However, in the general 3D 
case one cannot use the shear 
strain as the order parameter as we usually do  for the 1D case.  
One of the first approaches 
to deal with the 3D challenge was the Fr\'{e}mond  model. This model uses 
different expressions 
for free energy 
functions for different phases (see, for example, \cite{Fremond1990}).
All these expressions are essentially 
of the Landau-Ginsburg-type and contain the term 
${\gamma}/{2} \nabla {\rm tr} (\vec\epsilon)$ with 
$\gamma>0$ introduced in order to smooth  possibly very sharp 
spatial phase separation. 
In this paper we use a different approach proposed in \cite{Falk1990}.
This approach generalises  the classical Landau-Devonshire-Falk theory for 
 shape memory alloys to the 3D case. 
The free energy function, based on the expansion up to sixth order in  
 a single shear strain component \cite{Falk1980}, was extended to the three-dimensional case
using the group theoretical approach proposed in
\cite{Nittono1982}. 
In contrast to some other models (see, for example, \cite{Fremond1990})
strain-gradient terms are not involved in his expansion. 
We make use of this expansion and apply the following 
general representation of the free energy function
\begin{eqnarray}
\Psi({\vec\epsilon, \theta)= \psi^0(\theta) + 
\sum_{n=1}^{\infty} \psi^n 
(\vec\epsilon}, \theta ) \, ,
\label{6eq20}
\end{eqnarray}
where independent material parameters of the $n$-th order for $n=1,2,...$ 
are determined through strain invariants, ${\cal I}_j^n$ \, , as follows
\begin{eqnarray}
\psi^n = \sum_{j=1}^{j^n} \psi_j^n {\cal I}_j^n \quad \mbox{and} \quad 
\psi^0(\theta) = \psi_0(\theta) \, .
\label{6eq21}
\end{eqnarray}
The upper limit of the sum in (\ref{6eq21}), $j^n$, is the number 
of all invariant 
directions associated with 
a 
representation of the 48th order cubic symmetry group of the parent phase
(see details in \cite{Falk1990}).  
In order to adequately describe thermomechanical behaviour of 
shape-memory alloys we need to account for 6 terms in the sum of the 
expansion (\ref{6eq20}). In this case we have to determine  
 32 material parameters that 
make the application of formulae (\ref{6eq20})--(\ref{6eq21}) 
fairly complicated.
Using physically justified assumptions it is possible to reduce the number 
of required parameters. To achieve this, we make the same  assumptions 
as in \cite{Falk1990}. They 
 conclude  that odd degree invariants 
can be neglected in the expansion.
Taking invariants up to the sixth order results in a  
representation with only 10 material constants which may depend on 
temperature
\begin{eqnarray}
\Psi = \psi^0(\theta) + \sum_{j=1}^3 \psi_j^2 {\cal I}_j^2 + \sum_{j=1}^5 \psi_j^4 
{\cal I}_j^4
+ \sum_{j=1}^2 \psi_j^6 {\cal I}_j^6
\label{6eq24}
\end{eqnarray}
(we do not neglect the contribution of  $\psi^0(\theta)$ 
as  was done in \cite{Falk1990}).
The strain invariants ${\cal I}_i^n$ of second, forth and sixth orders of the 48th order
 cubic symmetry group  of the parent phase are 
\begin{eqnarray}
& &
\nonumber
{\cal I}_1^2= \frac{1}{9} ({\rm tr}(\epsilon_{ij}))^2, \; \; 
{\cal I}_2^2 = \frac{1}{12} 
(2\epsilon_{33} - \epsilon_{11} - \epsilon_{22})^2
+ \frac{1}{4} (\epsilon_{11} -\epsilon_{22})^2, \; \;
\\[10pt]
& &
\nonumber
{\cal I}_3^2 = 
\epsilon_{23}^2 + \epsilon_{13}^2 + \epsilon_{12}^2,
\; \;
{\cal I}_1^4 =({\cal I}_2^2)^2, \; \; {\cal I}_2^4 =\epsilon_{23}^4 + 
\epsilon_{13}^4 + \epsilon_{12}^4, \; \; {\cal I}_1^6=({\cal I}_2^2)^3
\\[10pt]
& &
\nonumber
{\cal I}_3^4 = ({\cal I}_3^2)^2, \; \; {\cal I}_4^4 = {\cal I}_2^2 
{\cal I}_3^2, \; \;
{\cal I}_5^4  = \epsilon_{23}^2 \left [\frac{1}{6} (2 \epsilon_{33} 
- \epsilon_{11} - \epsilon_{22})-
\frac{1}{2} (\epsilon_{11}- \epsilon_{22}) \right ]^2+ 
\\[10pt]
& &
\epsilon_{13}^2
\left [ \frac{1}{6} (2 \epsilon_{33} - \epsilon_{11} -
\epsilon_{22}) +
\frac{1}{2} (\epsilon_{11}- \epsilon_{22}) \right ]^2 + 
\frac{1}{9} \epsilon_{12}^2 (2 \epsilon_{33} - \epsilon_{11} 
- \epsilon_{22})^2,
\; \; 
\\[10pt]
& &
\nonumber
 {\cal I}_2^6 =  
\frac{1}{36} (2 \epsilon_{33} - \epsilon_{11} -
\epsilon_{22})^2 \left (
\frac{1}{36} (2 \epsilon_{33} - \epsilon_{11} -
\epsilon_{22})^2 - \frac{1}{4} ( \epsilon_{11}- \epsilon_{22})^2
\right ) ^2.
\end{eqnarray}
The ten material constants $\psi_j^n$ in (\ref{6eq24}) differ from alloy to alloy and we use 
coefficients 
derived for Cu-based alloys \cite{Falk1990} (units used here are 
consistent with those used in Section 7 for our  numerical results
 on the dynamics of shape-memory alloys):
\begin{eqnarray}
& &
\nonumber
\psi_1^2 =5.92 \times 10^6 {\rm \; g}/({\rm ms^2 cm}), \; \; \psi_2^2 = 
(1.41\times 10^5  + 46(\theta -300) ) {\rm \; g}/({\rm ms^2 cm}), \;
\\[10pt]
& &
\nonumber
\psi_3^2=(1.48 \times 10^6 - 940(\theta -300 ) ) {\rm \; g}/({\rm ms^2 cm}),
 \; \; \psi^0 = - \alpha_1 \theta \ln [(\theta - \theta_0)/\theta_0]
{\rm \; g}/({\rm ms^2 cm}), 
\\[10pt]
& &
\nonumber
\psi_1^4 = (-1.182 \times 10^8  + 3.55 \times 10^5 (\theta -300 ) )
{\rm \; g}/({\rm ms^2 cm}), \; \;
\\[10pt]
& &
\psi_2^4 =3.13 \times 10^9 {\rm \; g}/({\rm ms^2 cm}), \; \; 
\psi_3^4 =1.64 \times 10^9 {\rm \; g}/({\rm ms^2 cm}),
\\[10pt]
& &
\nonumber
\psi_4^4 =-5.53 \times 10^8 {\rm \; g}/({\rm ms^2 cm}), \; \; 
\psi_5^4 = -4.27 \times 10^8 {\rm \; g}/({\rm ms^2 cm}), 
\\[10pt]
& &
\nonumber
\psi_1^6 =3.35 \times 10^{10} {\rm \; g}/({\rm ms^2 cm}), \; \; 
\psi_2^6 = 3.71 \times 10^{11} {\rm \; g}/({\rm ms^2 cm}). 
\end{eqnarray}
Other material parameters are taken to be the same as those given in 
Section 7.
We are interested in the construction of an 
adequate model for the description of thermomechanical behaviour 
of thin slabs in
shape memory alloy materials.   
Starting from the 3D Falk-Konopka model and 
using centre manifold techniques (see, for example, \cite{Roberts1993})
we  derive systematically an 
accurate low-dimensional model for the dynamics of the slab.
The shape memory alloy is assumed to be of very large extent in the
$x=x_1$ direction compared to its thickness of $2b$ in the
$y=x_2$ direction ($-b<y<b$).
For the sake of convenience we use a new temperature variable
 $\theta'=\theta-\theta_0$ where here
$\theta_0=300$.
For simplicity of the analysis  we  assume  zero dissipation,
$\alpha=\mu=0$, and that there is no motion nor
dependence in the $x_3$ direction.

A model for the dynamics of modes which vary
slowly along the slab
is derived for the unforced dynamics, ${\bf F} = {\bf 0}$, $G=0$,
and when ``zero-stress \& thermal-insulation'' boundary
conditions are specified  on $y=\pm b$.
The derivation of boundary conditions in the ``long'' direction $x$ 
requires a quite delicate analysis \cite{Roberts1992} and these issues 
will not be address here. We only note  that 
``pinned \& insulating ends'' boundary conditions
may be used as a leading approximation.  
Modelling the long-wavelength, small-wavenumber modes along the slab,
we neglect all longitudinal variations
 and look for eigenvalues of the cross-slab
modes. It can be shown that generally there is a zero eigenvalue of 
multiplicity
five and all the rest are pure imaginary (as dissipation has been
omitted).
Thus there exists a sub-centre manifold based upon
these five modes (see \cite{Sijbrand85} for an existence 
theorem), called a slow manifold as these five modes evolve
slowly.
Note that being on a sub-centre manifold the models we construct only
have a weak assurance of asymptotic completeness, see the discussion
in \cite{Cox93b}.
The zero eigenvalue of multiplicity five corresponds to longitudinal
waves, large-scale bending, and one heat mode.
The leading order structure of the critical eigenmodes are constant
across the slab.
Thus letting an overbar denote the $y$ average,  the
amplitudes of the critical modes are chosen in the form
\begin{equation}
	U_i(x,t)=\overline{u_i}\,,\quad
	V_i(x,t)=\overline{v_i}\,,\quad
	\Theta'(x,t)=\overline{\theta'}\,.
	\label{Eampls}
\end{equation}
The low-dimensional model below is written in terms of these parameters.

The construction of the low-dimensional model is based upon the ansatz
that there exists a low-dimensional invariant manifold upon which the
amplitudes evolve slowly:
\begin{eqnarray}
&&
	u_i={\cal U}_i(\vec U,\vec V,\Theta')\,,\quad
	v_i={\cal V}_i(\vec U,\vec V,\Theta')\,,\quad
	\theta={\cal T}(\vec U,\vec V,\Theta')\,,
	\label{Eaman}
\\\mbox{where}&&
    \D t{U_i}=V_i\,,\quad
    \D t{V_i}=g_i(\vec U,\vec V,\Theta')\,,\quad
    \D t{\Theta'}=g_\theta(\vec U,\vec V,\Theta')\,.
    \label{Eamod}
\end{eqnarray}
This ansatz is substituted into the differential-algebraic
equations of 3D thermo-visco-elasticity (\ref{2eq12}), (\ref{2eq18}) 
and is solved to some
order in the small parameters $\partial_x$, $E=\|\vec U_x\|+\|\vec
V_x\|$ and $\vth=\|\Theta'\|$ with the computer algebra package 
{\textsc reduce}.
Thus, here we treat the strains as small, as measured by $E$, while
permitting asymptotically large displacements and velocities.
The displacement and temperature fields of the slow manifold, in terms
of the amplitudes and the scaled transverse coordinate $Y=y/b$, are
approximately
\begin{eqnarray}
	u_1 & \approx & U_1 -YbU_{2x} +0.15(3\,Y^2-1)b^2U_{1xx}\,,
	\label{Eslow1}  \\
	u_2 & \approx & U_2 -(0.9-3.05\E{-5}\Theta')YbU_{1x}
	+0.15(3\,Y^2-1)b^2U_{2xx}
	\nonumber\\&&{}
	-141\,Yb{U_{1x}}^3
	+1.00\E{-4}(3Y-Y^3)b^3{V_{1x}}^2U_{1x}\,,
	\label{Eslow2}  \\
	\theta & \approx & 300+\Theta'
	-2.43\E6(3Y-Y^3)b^3\left( V_{1x}U_{2xx} +U_{1x}V_{2xx} \right)
	\nonumber\\&&{}
	-25.1(7-30\,Y^2+15\,Y^4){V_{1x}}^3U_{1x}
	\,.
	\label{Eslowt}
\end{eqnarray}
These expressions have errors
$\Ord{E^5+\dx^{5/2}+\vth^{5/2}}$ where the notation
$\Ord{E^p+\dx^q+\vth^r}$ is used to denote terms involving $\dx^bE^a\vth^c$
such that $a/p+b/q+c/r\geq 1$. 
The mechanical and thermal field approximations represented by 
(\ref{Eslow1})--(\ref{Eslowt})
have   cross-slab structure.
In particular, the sideways deformation $u_2$ (which is a nonlinear function
of the longitudinal strains) of the shape memory alloy 
 feed back at
higher order to contribute to and 
complicate the longitudinal and thermal dynamics.

The model for the longitudinal dynamics on this slow manifold is
\begin{eqnarray}
	\rho\D t{V_1} & = & 2.97\E6 U_{1xx}+8.03\E5 b^2U_{1xxxx}
	\nonumber  
\\
& &
 {}+\dx\left[ (922\,\Theta'-0.0145{\Theta'}^2)U_{1x}
	 -(4.28\E9 -1.31\E7\Theta' ) {U_{1x}}^3
	 +7.12\E{11}{U_{1x}}^5
	\right.
\nonumber
\\
&&
\left.
{}	+(2820  -8.80\,\Theta')b^2 {V_{1x}}^2U_{1x}
	\right.
\left.
	+1.24\,b^4 {V_{1x}}^4U_{1x}
	-5.42\E4 b^2{V_{1x}}^2{U_{1x}}^3
	\right]
	\nonumber
\\
&&
{}+\Ord{E^8+\dx^4+\vth^4}\,. \label{Elongs}
\end{eqnarray}
The first line in the right-hand side of (\ref{Elongs}) describes 
linear dispersive elastic waves along
the slab, whereas the
 second line gives the temperature dependent quintic stress-strain
relation of the shape memory alloy.
Since $V_{1x}=U_{1xt}$, the remaining lines show effects upon this
stress-strain relation due to rates of change of the strain.

Note that to this order of truncation there is no coupling to the
bending modes of the slab which to the same error is simply the beam
equation
\begin{eqnarray}
	\rho\D t{V_2} & = & -9.91\E5 b^2U_{2xxxx}
	+\Ord{E^8+\dx^4+\vth^4}\,.
	 \label{Ebends}
\end{eqnarray}
There exists nonlinear coupling between the modes at higher order.

The corresponding energy equation for the temperature is
\begin{eqnarray}
	C_v\D t{\Theta'} & = & {\kappa}\Theta'_{xx}
+(2.77\E5 +914\,\Theta' -9.25\,{\Theta'}^2)U_{1x}V_{1x}
	\nonumber\\&&{}
	 +(3.94\E9 +1.26\E7\,\Theta'  ) V_{1x}{U_{1x}}^3
	 -(57.3  +0.0117\,\Theta' )b^2 {V_{1x}}^3 U_{1x}
	\nonumber\\&&{}
	 +1.68\E{12} V_{1x}{U_{1x}}^5
	 -1.58\E6 b^2 {V_{1x}}^3{U_{1x}}^3
	 -0.0203\,b^4 {V_{1x}}^5U_{1x}
	\nonumber\\&&{}
	 +1.63\E4 b^2 U_{1xx}V_{1xx}
	 +9.22\E4 b^2 U_{2xx}V_{2xx}
+\dx^2\left[ -8151\,b^2 U_{1x}V_{1x}\right]
	\nonumber\\&&{}+\Ord{E^8+\dx^4+\vth^4}\,. \label{Ethes}
\end{eqnarray}
The first line in (\ref{Ethes}) describes the diffusion of heat 
generated or absorbed by mechanical strains, ${\Theta}U_{1x}V_{1x}$.
However, in the thin slab the internal pattern of strains causes a much more complicated distribution of heating and cooling as summarised by the remaining
lines.
 It is expected that virtually all of these should be
retained in order to be consistent with the quintic stress-strain of
the longitudinal wave equation.
Computational experiments with the model derived 
in this section will be presented elsewhere.

\section{Acknowledgements}
The authors were supported by ARC
Small Grant 
179406. 
Special thanks go to Kerryn Thomas
 who  contributed  the results of 
computational experiments.

%
%

%
%
\paperend
\end{document}